\documentclass[12pt,oneside,reqno]{amsart}
\usepackage{}
\usepackage{amssymb}
\usepackage{bbm}
\usepackage{cases}
\usepackage{amsmath}
\usepackage{graphicx}
\usepackage{mathrsfs}
\usepackage{stmaryrd}
\usepackage{amsfonts}
\usepackage{enumerate,amsmath,amssymb,amsthm}

\numberwithin{equation}{section}

\newcommand{\be}{\begin{eqnarray}}
\newcommand{\ee}{\end{eqnarray}}
\newcommand{\ce}{\begin{eqnarray*}}
\newcommand{\de}{\end{eqnarray*}}
\newtheorem{theorem}{Theorem}[section]
\newtheorem{lemma}[theorem]{Lemma}
\newtheorem{remark}[theorem]{Remark}
\newtheorem{definition}[theorem]{Definition}
\newtheorem{proposition}[theorem]{Proposition}
\newtheorem{Examples}[theorem]{Example}
\newtheorem{corollary}[theorem]{Corollary}

\def\e{{\mathrm{e}}}

\def\[{{\Big[}}
\def\]{{\Big]}}
\def\<{{\langle}}
\def\>{{\rangle}}
\def\({{\Big(}}
\def\){{\Big)}}

\def\bx{{\mathbf{x}}}

\def\dif{{\mathord{{\rm d}}}}

\def\no{\nonumber}
\def\={&\!\!=\!\!&}
\def\bt{\begin{theorem}}
\def\et{\end{theorem}}
\def\bl{\begin{lemma}}
\def\el{\end{lemma}}
\def\br{\begin{remark}}
\def\er{\end{remark}}

\def\bd{\begin{definition}}
\def\ed{\end{definition}}
\def\bp{\begin{proposition}}
\def\ep{\end{proposition}}
\def\bc{\begin{corollary}}
\def\ec{\end{corollary}}
\def\bx{\begin{Examples}}
\def\ex{\end{Examples}}

\def\cI{{\mathcal I}}
\def\cJ{{\mathcal J}}

\def\cQ{{\mathcal Q}}

\def\cV{{\mathcal V}}
\def\cW{{\mathcal W}}

\def\mD{{\mathbb D}}

\def\mK{{\mathbb K}}

\def\mN{{\mathbb N}}

\def\mR{{\mathbb R}}

\def\sG{{\mathscr G}}

\def\sL{{\mathscr L}}

\def\sQ{{\mathscr Q}}

\def\geq{\geqslant}
\def\leq{\leqslant}

\allowdisplaybreaks

\begin{document}

\title{Sharp heat kernel estimates for spectral fractional Laplacian perturbed by gradient}

\date{}

\author{Renming Song,\ \ Longjie Xie \ \  and \ \ Yingchao Xie}

\address{Renming Song:
Department of Mathematics, University of Illinois,
Urbana, IL 61801, USA\\
Email: rsong@illinois.edu
 }

\address{Longjie Xie:
School of Mathematics and Statistics, Jiangsu Normal University,
Xuzhou, Jiangsu 221000, P.R.China\\
Email: xlj.98@whu.edu.cn
 }
\address{Yingchao Xie:
School of Mathematics and Statistics, Jiangsu Normal University,
Xuzhou, Jiangsu 221000, P.R.China\\
Email: ycxie@jsnu.edu.cn
 }

\thanks{Research of R. Song is supported by the Simons Foundation (\#429343, Renming Song). L. Xie is supported by NNSF of China (No. 11701233) and NSF of Jiangsu (No. BK20170226). Y. Xie is supported by NNSF of China (No. 11771187). The Project Funded by the PAPD of Jiangsu Higher Education Institutions is also gratefully acknowledged}

\begin{abstract}
By using Duhamel's formula, we prove sharp two-sided estimates for the heat kernel of spectral fractional Laplacian with time-dependent gradient perturbation in
bounded $C^{1,1}$ domains.
Moreover, we also obtain gradient estimate as well as  H\"older continuity of the gradient of the heat kernel.
\bigskip

\end{abstract}

\maketitle \rm

\section{Introduction and main results}
Let $W_{t}$ be a Brownian motion in $\mathbb{R}^{d}$ $(d\geq1)$ with generator $\Delta$
and $T_{t}$ be an independent $\alpha/2$-stable subordinator
with $\alpha\in(0,2)$. Then the subordinate process
$X_{t}:=W_{T_{t}}$ is an isotropic $\alpha$-stable
process and its infinitesimal generator is the fractional Laplacian operator $-(-\Delta^{\alpha/2})$ which is given by
$$
-(-\Delta^{\alpha/2})f(x):= \int_{\mR^d}\big[f(x+z)-f(x)-1_{|z|\leq
1}z\cdot\nabla f(x)\big]\frac{c_{d,\alpha}}{|z|^{d+\alpha}}\dif z,\quad\forall f\in C^2_c(\mR^d),
$$
where $c_{d,\alpha}$ is a positive constant.
It is well known that the heat kernel $p(t,x,y)$ of $-(-\Delta^{\alpha/2})$ (which is
also the transition density of $X:=(X_t)_{t\geq 0}$)
has the following estimates: for every $t>0$ and $x,y\in\mR^d$,
\begin{align}
p(t,x,y)\asymp\left(t^{-d/\alpha}\wedge\frac{t}{|x-y|^{d+\alpha}}\right).  \label{Heat}
\end{align}
Here and below, for two non-negative functions $f$ and $g$, the notation $f \asymp g$ means that there are positive constants $c_1$ and $c_2$ such that $c_1g(x)\leq f(x)\leq c_2g(x)$ in the common domain of $f$ and $g$.

In \cite{Bo-Ja0}, by using Duhamel's formula,
Bogdan and Jakubowski studied the following perturbation of
$-(-\Delta^{\alpha/2})$ by a gradient operator:
$$
\sL^b:=-(-\Delta^{\alpha/2})+b(x)\cdot\nabla,\quad\alpha\in(1,2),
$$
where $b=(b^1, \cdots, b^d): \mR^d\to\mR^d$ with $b^j$, $j=1, \dots, d$, belonging
to the Kato class ${\bf K}^{\alpha-1}_d$ defined
as follows: for $\gamma>0$,
$$
{\bf K}^{\gamma}_d:=\left\{f\in L^1_{loc}(\mR^d): \lim_{r\downarrow 0}\sup_{x\in\mR^d}
\int_{B(x,r)}\frac{|f(y)|}{|x-y|^{d-\gamma}}\dif y=0\right\},
$$
and $B(x,r)$ denotes the open ball centered at $x\in\mathbb{R}^{d}$ with radius $r$.
Let $p^b(t, x, y)$ be the heat kernel of $\sL^b$.
Small time sharp two-sided estimates for $p^b$,  of the form \eqref{Heat},
were established in \cite[Theorems 1 and 2]{Bo-Ja0}.
The key points of the perturbation method used in \cite{Bo-Ja0} are that, on one hand, a nice bound on $\nabla_x p(t,x,y)$ is known, and on the other hand, the following 3-P inequality concerning $p(t,x,y)$ holds: there exists $C_0>0$ such that for any $0<s<t$ and $x,y,z\in \mR^d$,
\begin{align}
\frac{p(t-s,x,z)p(s,z,y)}{p(t,x,y)}\leq C_0\Big(p(t-s,x,z)+p(s,z,y)\Big).  \label{ppp}
\end{align}
See also \cite{C-K-K,Ch-Ku2,C-K-W,J-S,K-S,W-Z,Xi-Zh} and the references therein for two-sided heat kernel estimates of more general non-local operators in the whole space $\mR^d$.

\vspace{2mm}
Let $D$ be an open subset of $\mR^d$, we can kill the process $X$ upon exiting D and obtain a
subprocess $X^D$ known as the killed isotropic $\alpha$-stable process.
The infinitesimal generator of $X^D$ is the Dirichlet fractional Laplacian $-(-\Delta)^{\alpha/2}|_D$,
that is, the fractional Laplacian with zero exterior condition.
Due to the complication near the boundary, two-sided
estimates for the Dirichlet heat kernel of $-(-\Delta)^{\alpha/2}|_D$ (or equivalently, the transition density of $X^D$) are much more difficult to obtain. To state the related results, we first recall that an open set $D$ in $\mR^{d}$ is said to be $C^{1,1}$ if there exist $r_{0}>0$ and
$\Lambda>0$ such that for every $Q\in\partial D$, there exist a $C^{1,1}$-function $\phi=\phi_{Q}: \mathbb{R}^{d-1}\rightarrow\mathbb{R}$ satisfying
$\phi(0)=\nabla\phi(0)=0$, $\|\nabla\phi\|_{\infty}\leq\Lambda$, $|\nabla\phi(x)-\nabla\phi(z)|\leq\Lambda|x-z|$ and an orthonormal coordinate system $y= (y_{1}, \cdots, y_{d-1},
y_{d}):=(\tilde{y}, y_{d})$ such that $B(Q, r_{0})\cap D=B(Q,r_{0})\cap \{y: y_{d}>\phi(\tilde{y})\}$. The pair $(r_{0},\Lambda)$ is called the characteristics of the $C^{1,1}$ open set $D$. In \cite{C-K-S-2}, Chen, Kim and Song proved that when $D$ is a $C^{1,1}$ open set in $\mR^d$, the heat kernel $p^D(t,x,y)$ of $-(-\Delta)^{\alpha/2}|_D$ has the following two-sided
estimates: for every $T>0$ and $(t,x,y)\in(0,T]\times D\times D$,
\begin{align}
p^D(t,x,y)\asymp
\left(1\wedge\frac{\rho(x)^{\alpha/2}}{\sqrt{t}}\right)\left(1\wedge\frac{\rho(y)^{\alpha/2}}{\sqrt{t}}\right)p(t,x,y),\label{Heat2}
\end{align}
where $\rho(x)$ denotes the distance between $x$ and $D^{c}$.

By reversing the order of subordination and killing, one can obtain a process $Y^D$ which is different from
$X^D$. More precisely, we first kill the Brownian motion $W$ at $\tau_{D}$, the first exit time of $W$ from $D$, and then subordinate the killed Brownian motion $W^D$ using the independent $\alpha/2$-stable subordinator $T_{t}$. That is, $Y^D:=(W^{D})_{T_{t}}$ is defined as
\begin{align*}
Y^D_{t}:=
\begin{cases}
W_{T_{t}}, & T_{t}<\tau_{D}\\
\partial, & T_{t}\geq\tau_{D}
\end{cases}
=
\begin{cases}
W_{T_{t}}, & t<A_{\tau_{D}}\\
\partial, & t\geq A_{\tau_{D}},
\end{cases}
\end{align*}
where $\partial$ is a cemetery state, $A_{t}:=\inf\{s>0: T_{s}\geq t\}$ is the inverse of $T$ and the last equality follows from the fact $\{T_{t}<\tau_{D}\}=\{t<A_{\tau_{D}}\}$.
The process $Y^D$ is called a subordinate killed Brownian motion. For the differences and relationship between
the processes $X^{D}$ and $Y^D$, see \cite{S-V2}. The infinitesimal generator of $Y^D$ is the spectral fractional Laplacian $-(-\Delta|_{D})^{\alpha/2}$, which is defined as a fractional power of the negative Dirichlet Laplacian. It is a very useful object in analysis and partial differential equations (see \cite{B-S-V,P-A,S-Z}) and has been intensively studied (see \cite{A-D,D-M-Z,G-P-R-S-S-V,S-V} and the references therein).
When $D$ is a bounded $C^{1,1}$ domain,
the following sharp estimates for the heat kernel $r^D(t,x,y)$ of $-(-\Delta|_{D})^{\alpha/2}$ (which is also the transition density of $Y^D$) were obtained in \cite[Theorem 4.7]{Song}: for every $T>0$ and $(t,x,y)\in(0,T]\times D\times D$,
\begin{align*}
r^D(t,x,y)\asymp\left(1\wedge\frac{\rho(x)\rho(y)}{\left(|x-y|+t^{1/\alpha}\right)^{2}}\right)p(t,x,y).
\end{align*}
In Lemma \ref{sharp} below, we will give the following
alternative form of the estimates above: for $(t,x,y)\in(0,T]\times D\times D$,
\begin{align*}
r^D(t,x,y)\asymp\left(1\wedge\frac{\rho(x)}{|x-y|+t^{1/\alpha}}\right)\left(1\wedge\frac{\rho(y)}{|x-y|+t^{1/\alpha}}\right)p(t,x,y),
\end{align*}
which is more convenient to use.

\vspace{2mm}
Gradient perturbations of Dirichlet operators have also been widely studied in recent years.
In \cite{C-K-S-1}, Chen, Kim and Song studied the following gradient perturbation of the Dirichlet fractional Laplacian:
$$
\sL^{b,D}:=\left(-(-\Delta)^{\alpha/2}+b(x)\cdot\nabla\right)|_{D},\quad\alpha\in(1,2).
$$
Under the condition that $b\in{\bf K}^{\alpha-1}_d$ and $D$ is a
bounded $C^{1,1}$ open set in $\mathbb{R}^{d}$ with $d\geq2$, Chen,
Kim and Song \cite[Theorem 1.3]{C-K-S-1} showed that the heat kernel $p^{b,D}(t,x,y)$ of $\sL^{b,D}$ has the same
estimates as in (\ref{Heat2}). This result was generalized to unbounded $C^{1,1}$ open sets by \cite{P-R}. Unlike the whole space case, there was no good estimate on
$\nabla_x p^D(t, x, y)$, thus \cite{C-K-S-1, P-R} used Duhamel's formula
for the Green function and the probabilistic road-map designed in \cite{C-K-S-2} for establishing
the estimates \eqref{Heat2}.

In the recent paper \cite{K-R},  Kulczycki and Ryznar  proved the following gradient estimate
for $p^D(t, x, y)$: for any $T>0$, there exists a constant $C=C(d,T)>0$ such that for any $(t, x, y)\in (0, T]\times D\times D$,
\begin{align*}
|\nabla_x p^{D}(t,x,y)|\leq  \frac{C}{\rho(x)\wedge
t^{1/\alpha}}p^{D}(t,x,y).
\end{align*}
Using this result, we gave, in the recent preprint \cite{C-S-X-X}, a direct proof of the main results in \cite{C-K-S-1, P-R} by using Duhamel's formula, with drift $b=(b^1, \cdots, b^d): D\to \mR^d$, where each $b^j$, $j=1, \dots, d$, belongs to the following  Kato class:
$$
{\bf K}_{D}^{\alpha-1}:=\left\{f\in L^1_{loc}(D): \lim_{r\downarrow0}\sup_{x\in D}\int_{D\cap
B(x,r)}\frac{|f(y)|}{|x-y|^{d+1-\alpha}} \dif y=0\right\}.
$$
Moreover, we also obtain a gradient estimate for $p^{b,D}(t,x,y)$. Notice that by H\"older's inequality, $L^p(D)\subseteq {\bf K}_{D}^{\alpha-1}$ provided $d/(\alpha-1)<p\leq \infty$.

\vspace{2mm}
The aim of this paper is to study the following perturbation of spectral fractional Laplacian by a time-dependent gradient operator:
$$
\sL^{D,b}:=-(-\Delta|_D)^\frac{\alpha}{2}+b(t,x)\cdot\nabla,\quad\alpha\in(1,2),
$$
with $b(t, x)=(b^1(t, x), \cdots, b^d(t, x)): (0, \infty)\times D\to\mR^d$ satisfying certain conditions.
We shall derive sharp two-sided estimates for the heat kernel $r^{D, b}(t, x, y)$
of $\sL^{D,b}$ in bounded $C^{1,1}$ domains.
Moreover, we also obtain a
gradient estimate as well as the H\"older continuity of the gradient of $r^{D, b}(t, x, y)$, which are of independent interest.

To state our main result, let us first introduce our local Kato
class of space-time functions used in this paper.

\bd\label{Definition}
Let $D$ be a domain in $\mR^{d}$
and $\gamma\geq 0$. For a real-valued function $f$ on $(0,
\infty)\times D$, we define for every $\delta>0$,
\begin{align*}
K^{\gamma}_{f}(\delta):=\sup_{t>0,x\in
D}\delta^{\gamma/\alpha}\!\int_{0}^{\delta}\!\!\!\int_{D}&\big[s^{-\gamma/\alpha}+(\delta-s)^{-\gamma/\alpha}\big]\left(1\wedge\frac{\rho(y)}{|x-y|+s^{1/\alpha}}\right)\\
&\quad\times\frac{s}{(|x-y|+s^{1/\alpha})^{d+\alpha+1}}\cdot|f(t\pm s,y)|\dif y\dif s.
\end{align*}
We say that the function $f$ belongs to thel Kato class $\mK_{D}^{\gamma}$ if
$\lim_{\delta\downarrow0}K_{f}^{\gamma}(\delta)=0$.
\ed

\br
We note that our Kato class is time-dependent, which is needed when consider parabolic problems, see \cite{J-S,ZhQ}. One can easily check that if $0\leq\gamma_1<\gamma_2$, then $\mK_D^{\gamma_2}\subseteq\mK_D^{\gamma_1}$.
By Lemma \ref{include} below, we have that ${\bf K}_D^{\alpha-1}\subset \mK_D^{0}$
and that, for  $1<p,q\leq \infty$, $L^q(\mR;L^p(D))\subseteq \mK_D^{\gamma}$ provided $\tfrac{d}{\alpha p}+\tfrac{1}{q}<1-\tfrac{1+\gamma}{\alpha}$.
\er

In the remainder of this paper, we always assume that $b=(b^1, \cdots, b^d): (0, \infty)\times D\to\mR^d$ and each $b^j$, $j=1, \dots, d$, belongs to $\mK_{D}^0$ .

By Duhamel's formula, the heat
kernel $r^{D,b}(s,x;t,y)$ of $\sL^{D,b}$ should satisfy the following integral equation: for $0\leq s<t$ and $x,y\in D$,
\begin{align}\label{Duhamel}
r^{D,b}(s,x;t,y)=r^D(t-s,x,y)+\int_s^t\!\!\!\int_{D}r^{D,b}(s,x;r,z)b(r,z)\cdot\nabla_z
r^D(t-r,z,y)\dif z\dif r,
\end{align}
or
\begin{align}\label{duhamel}
r^{D,b}(s,x;t,y)=r^{D}(t-s,x,y)+\int_s^t\!\!\!\int_D
r^{D}(r-s,x,z)b(r,z)\cdot\nabla_z r^{D,b}(r,z;t,y)\dif z\dif r.
\end{align}
Notice that in (\ref{Duhamel}) the derivative of the unknown heat kernel is not involved, and hence easier to solve.
While (\ref{duhamel}) is connected directly to the mild solutions of the corresponding parabolic equations, and from which one can easily derive the H\"older continuity of the gradient of the unknown heat kernel.
For convenience, we define for $t>0$ and $x,y\in D$,
\begin{equation}\label{qq}
q^D(t, x, y):=\left(1\wedge\frac{\rho(x)}{|x-y|+t^{1/\alpha}}\right)\left(1\wedge\frac{\rho(y)}{|x-y|+t^{1/\alpha}}\right)p(t,x,y).
\end{equation}
The following is the main result of this paper.

\bt\label{main}
Let $D$ be a bounded $C^{1,1}$ domain in $\mathbb{R}^{d}$
and $b\in\mK^{0}_{D}$.
Then there exists a unique function $r^{D,b}(s,x;t,y)$ on $(0,
\infty)\times D\times D$ satisfying (\ref{Duhamel}) such that:

\begin{enumerate}[(i)]

\item (Two-sided estimates)
for any $\delta>0$,
there exists a constant $C_1>1$ such that for all $0\leq s<t\leq s+\delta$
and $x, y\in D$, we have
\begin{align}\label{estimate}
C_1^{-1}q^D(t-s,x,y)\leq r^{D,b}(s,x;t,y)\leq C_1 q^D(t-s,x,y);
\end{align}

\item (Gradient estimate)
for any $\delta>0$, there exists a constant $C_2>0$ such that for all $0\leq s<t\leq s+\delta$
and $x, y\in D$,
\begin{align}\label{vi}
|\nabla_x r^{D,b}(s,x;t,y)|\leq C_2
\frac{1}{\rho(x)\wedge\left(|x-y|+(t-s)^{1/\alpha}\right)}q^D(t-s,x,y),
\end{align}
and $r^{D,b}(s,x;t,y)$ also satisfies (\ref{duhamel});

\item (C-K equation) for all $0\leq s<r<t$ and $x,y\in D$, the following Chapman-Kolmogorov's equation holds:
\begin{align}\label{eq21}
\int_{D}r^{D,b}(s,x;r,z)r^{D,b}(r,z;t,y)\dif z=r^{D,b}(s,x;t,y);
\end{align}

\item (Generator) for any $f\in C^2_c(D)$, we have
\begin{align}\label{eq23}
R_{s,t}^{D,b}f(x)=f(x)+\int_s^tR_{s,r}^{D,b}\mathcal{L}^{D,b}f(x)\dif r,
\end{align}
where $R^{D,b}_{s,t}f(x):=\int_{D}r^{D,b}(s,x;t,y)f(y)\dif y$;

\item (Continuity) for any uniformly continuous function $f(x)$ with compact supports, we have
\begin{align}
\lim_{t\downarrow s} \|R^{D,b}_{s,t}f-f\|_\infty=0;   \label{con}
\end{align}

\item (H\"older continuity)
if we further assume that $b\in \mK_D^{\gamma}$ for some $\gamma\in (0, \alpha-1)$,
then for any $\delta>0$,
there exists a constant $C_3>0$ such that for any $0\leq s<t\leq s+\delta$
and $x,x',y\in D$, we have
\begin{align}
|\nabla_x r^{D,b}& (s,x;t,y)-\nabla_x r^{D,b}(s,x';t,y)| \leq C_3|x-x'|^\gamma (t-s)^{-\gamma/\alpha}\no\\
&\times\frac{1}{\rho(\widetilde x)\wedge\left(|\widetilde x-y|+(t-s)^{1/\alpha}\right)}q^D(t-s,\widetilde x,y),  \label{ho}
\end{align}
where $\widetilde x$ stands the point among $x$ and $x'$ which is closer to $y$.
\end{enumerate}
\et

We remark that the gradient estimates (\ref{vi}) and (\ref{ho}) are new even in the case $b\equiv0$.
We now briefly describe the main idea of our argument.
Due to the difference between the processes $Y^D$ and $X^D$,
the method used in \cite{C-K-S-1,P-R}  does not work for $\sL^{D,b}$.
Instead, we will use Duhamel's formula (\ref{Duhamel}) to obtain the sharp two-sided estimates of the heat kernel. As mentioned before, two main ingredients are needed: the gradient estimate for $r^D(t,x,y)$ and the corresponding 3-P inequality, both of which are unknown. In fact, by Remark \ref{pp} below, we shall see that the 3-P inequality of the form (\ref{ppp}) does not hold for the heat kernel $r^D(t,x,y)$. Because of these, we will first derive an estimate on $\nabla_x r^D(t,x,y)$, and then establish a generalized 3-P type inequality for $r^D(t,x,y)$.
The gradient estimate and the H\"older estimate for $r^{D,b}(s,x;t,y)$ follow as easy by-products of our perturbation argument.

\vspace{2mm}
The remainder of this paper is organized as follows. In Section 2, we
prepare some important inequalities for $r^D(t,x,y)$ and derive its first and second order gradient estimates. The proof of the
main result, Theorem \ref{main}, will be given in Section 3.

We conclude this introduction by spelling out some conventions that
will be used throughout this paper. The letter C with or without
subscripts will denote an unimportant constant and $f\preceq g$
means that $f\leq Cg$ for some $C\geq1$. The letter $\mathbb{N}$
will denote the collection of positive integers, and
$\mathbb{N}_{0}=\mathbb{N}\cup\{0\}$. We will use $:=$ to denote a
definition, and we assume that all the functions
considered in this paper are Borel measurable.

\section{Estimates for $r^{D}(t,x,y)$}

In the remainder of this paper, $D$ denotes a bounded $C^{1,1}$ domain in $\mR^d$.
For simplicity, we first introduce some functions for latter
use. Given $d\geq 1$, $\vartheta\in\mR$ and $\alpha\in(0,2]$, we define for $t>0$ and $x,y\in D$,
$$
\varrho_{d}^\vartheta(t,x):=\frac{t^\vartheta}{(|x|+t^{1/\alpha})^{d+\alpha}},
$$
and
\begin{align}
\hat q_\alpha(t,x,y):=1\wedge\frac{\rho(x)}{|x-y|+t^{1/\alpha}},\quad q_\alpha(t,x,y):=\hat q_\alpha(t,x,y)\hat q_\alpha(t,y,x).\label{q}
\end{align}
Then we have $p(t, x, y)\asymp \varrho_{d}^1(t,x-y)$ and $q^D(t, x, y)=q_\alpha(t,x,y)
p(t, x, y)$.

We will first establish a generalized 3-P type inequality for $r^D(t,x,y)$, and then derive its first and second order
gradient estimates, which will be essential in
constructing the solution to the integral equation (\ref{Duhamel}).

\subsection{Generalized 3-P inequality}
Let $T>0$ be fixed. Recall that $r^D(t,x,y)$ is the heat kernel of
$-(-\Delta|_D)^{\frac{\alpha}{2}}$, and for any $t\in(0, T]$ and $x, y\in
D$, we have
\begin{align*}
r^D(t,x,y)\asymp\left(1\wedge\frac{\rho(x)\rho(y)}{(|x-y|+t^{1/\alpha})^{2}}\right)\varrho^1_{d}(t,x-y).
\end{align*}
The estimates above are not very convenient for our application since $\rho(x)$ and $\rho(y)$
are intertwined together.
We prove the following result.

\bl\label{sharp}
For any $t\in(0,
T]$ and $x, y\in D$, we have
\begin{align}
r^D(t,x,y)\asymp q_\alpha(t,x,y)\varrho^1_{d}(t,x-y)\asymp q^D(t, x, y).
\label{er}
\end{align}
\el
\begin{proof}
The second comparison follows immediately the sentence after \eqref{q}. So we will only prove
the first comparison.
It is obvious that
$$
q_\alpha(t,x,y)\preceq 1\wedge\frac{\rho(x)\rho(y)}{(|x-y|+t^{1/\alpha})^{2}}.
$$
Thus we only need to show that
\begin{align}
q_\alpha(t,x,y)\succeq 1\wedge\frac{\rho(x)\rho(y)}{(|x-y|+t^{1/\alpha})^{2}}.    \label{00}
\end{align}
One can easily see that the above inequality holds when
$$
\rho(x)\vee\rho(y)\leq |x-y|+t^{1/\alpha}\quad\text{or}\quad\rho(x)\wedge\rho(y)\geq|x-y|+t^{1/\alpha}.
$$
By symmetry, it suffices to prove (\ref{00}) in the case when
$$
\rho(x)\leq |x-y|+t^{1/\alpha}\leq\rho(y).
$$
Using the fact that $\rho(y)\leq\rho(x)+|x-y|$, we can deduce
\begin{align*}
1\wedge\frac{\rho(x)\rho(y)}{(|x-y|+t^{1/\alpha})^{2}}&\leq 1\wedge\frac{\rho(x)(\rho(x)+|x-y|)}{(|x-y|+t^{1/\alpha})^{2}}\\
&\leq 1\wedge\frac{\rho(x)^2}{(|x-y|+t^{1/\alpha})^{2}}+1\wedge\frac{\rho(x)\cdot|x-y|}{(|x-y|+t^{1/\alpha})^{2}}\\
&\preceq 1\wedge\frac{\rho(x)}{|x-y|+t^{1/\alpha}},
\end{align*}
which implies the desired result.
\end{proof}

\br\label{pp}
By (\ref{er}) and the same argument as in \cite[Remark 2.3]{C-K-S-3}, one can see that for all $t/4<s<3t/4$ and $x,y,z\in D$ with $2|x-y|\geq |x-z|+|z-y|$, it holds that
\begin{align*}
\frac{r^D(t+s,x,y)[r^D(t,x,z)+r^D(s,z,y)]}{r^D(t,x,z)r^D(s,z,y)}&\preceq \left(\frac{\rho(x)[\rho(z)+|x-y|+(t+s)^{1/\alpha}]}{\rho(z)[\rho(x)+|x-y|+(t+s)^{1/\alpha}]}\right)\\
&\quad+\left(\frac{\rho(y)[\rho(z)+|x-y|+(t+s)^{1/\alpha}]}{\rho(z)[\rho(y)+|x-y|+(t+s)^{1/\alpha}]}\right),
\end{align*}
which goes to zero as $\rho(x)=\rho(y)\rightarrow0$. This means that, unlike (\ref{ppp}), the inequality
$$
\frac{r^D(t,x,z)r^D(s,z,y)}{r^D(t+s,x,y)}\preceq r^D(t,x,z)+r^D(s,z,y)
$$
can not be true for all $t,s>0$ and $x,y,z\in D$, even for balls.
\er

We now proceed to prove a generalized 3-P type inequality for $r^D(t,x,y)$. Let us start with the following result.
\bl For any $t,s\geq 0$ and
$x,y,z\in D$, we have
\begin{align}
\frac{q_\alpha(t,x,z)q_\alpha(s,z,y)}{q_\alpha(t+s,x,y)}\preceq [\hat q_\alpha(t,z,x)]^2+[\hat q_\alpha(s,z,y)]^2. \label{3q}
\end{align}

\el
\begin{proof}
Note that, for any $a, b>0$, it holds that
\begin{equation}\label{elmentary}
1\wedge \frac{a}{b}\asymp \frac{a}{a+b}.
\end{equation}
Thus,
\begin{align*}
\frac{\hat q_\alpha(t,x,z)\hat q_\alpha(s,y,z)}{q_\alpha(t+s,x,y)}&\asymp \frac{\left((t+s)^{1/\alpha}+|x-y|+\rho(x)\right)
\left((t+s)^{1/\alpha}+|x-y|+\rho(y)\right)}{\left(t^{1/\alpha}+|x-z|+\rho(x)\right) \left(s^{1/\alpha}+|z-y|+\rho(y)\right)}\\
&\preceq 1+\frac{t^{1/\alpha}+|x-z|}{s^{1/\alpha}+|z-y|+\rho(y)}+\frac{s^{1/\alpha}+|z-y|}{t^{1/\alpha}+|x-z|+\rho(x)}.
\end{align*}
By  (\ref{q}), we have
\begin{align*}
\cI:=\frac{q_\alpha(t,x,z)q_\alpha(s,z,y)}{q_\alpha(t+s,x,y)}&=\frac{\hat q_\alpha(t,x,z)\hat q_\alpha(s,y,z)}{q_\alpha(t+s,x,y)}\hat q_\alpha(t,z,x)\hat q_\alpha(s,z,y)\\
&\preceq \hat q_\alpha(t,z,x)\hat q_\alpha(s,z,y)+\frac{\rho(z)}{s^{1/\alpha}+|z-y|+\rho(y)}\hat q_\alpha(s,z,y)\\
&\qquad+\frac{\rho(z)}{t^{1/\alpha}+|x-z|+\rho(x)}\hat q_\alpha(t,z,x).
\end{align*}
Using the fact
$$
\rho(x)+|x-z|\asymp\rho(z)+|x-z|,
$$
we further get that
$$
\frac{\rho(z)}{t^{1/\alpha}+|x-z|+\rho(x)}\asymp\frac{\rho(z)}{t^{1/\alpha}+|x-z|+\rho(z)}\asymp \hat q_\alpha(t,z,x),
$$
and similarly
$$
\frac{\rho(z)}{s^{1/\alpha}+|z-y|+\rho(y)}\asymp\frac{\rho(z)}{s^{1/\alpha}+|z-y|+\rho(z)}\asymp \hat q_\alpha(s,z,y).
$$
Thus, we have
\begin{align*}
\cI\preceq \hat q_\alpha(t,z,x)\hat q_\alpha(s,z,y)+[\hat q_\alpha(t,z,x)]^2+[\hat q_\alpha(s,z,y)]^2\preceq [\hat q_\alpha(t,z,x)]^2+[\hat q_\alpha(s,z,y)]^2.
\end{align*}
The proof is finished.
\end{proof}

As a direct consequence, we can obtain the following generalized 3-P type inequality for $r^D(t,x,y)$.

\bl
Let $T>0$. For any $0\leq s,t\leq T $ and $x, y, z\in D$, it holds that
\begin{align}
\frac{r^D(t,x,z)r^D(s,z,y)}{r^D(t+s,x,y)}\preceq (t\wedge s)\Big(&[\hat q_\alpha(t,z,x)]^2\varrho_{d}^0(t,x-z)+[\hat q_\alpha(s,z,y)]^2\varrho_{d}^0(s,z-y)\Big).   \label{3r}
\end{align}
\el
\begin{proof}
Combining (\ref{er}) and (\ref{3q}), we get that
\begin{align*}
\cJ:=\frac{r^D(t,x,z)r^D(s,z,y)}{r^D(t+s,x,y)}\preceq \Big([\hat q_\alpha(t,z,x)]^2+[\hat q_\alpha(s,z,y)]^2\Big)
\frac{\varrho_{d}^1(t,x-z)\varrho_{d}^1(s,z-y)}{\varrho_{d}^1(t+s,x-y)}.
\end{align*}
Note that
\begin{align*}
\big(|x-y|+(t+s)^{1/\alpha}\big)^{d+\alpha}\preceq \big(|x-z|+t^{1/\alpha}\big)^{d+\alpha}+\big(|z-y|+s^{1/\alpha}\big)^{d+\alpha}.
\end{align*}
Thus
\begin{align}
\frac{\varrho_{d}^1(t,x-z)\varrho_{d}^1(s,z-y)}{\varrho_{d}^1(t+s,x-y)}&=\frac{t\cdot s}{t+s}\cdot\frac{\varrho_{d}^0(t,x-z)\varrho_{d}^0(s,z-y)}{\varrho_{d}^0(t+s,x-y)}\no\\
&\preceq (t\wedge
s)\Big(\varrho_{d}^0(t,x-z)+\varrho_{d}^0(s,z-y)\Big).  \label{333}
\end{align}
Hence
\begin{align*}
\cJ&\preceq (t\wedge s)\Big([\hat q_\alpha(t,z,x)]^2+[\hat q_\alpha(s,z,y)]^2\Big)\Big(\varrho_{d}^0(t,x-z)+\varrho_{d}^0(s,z-y)\Big)\\
&\preceq (t\wedge s)\Big([\hat q_\alpha(t,z,x)]^2\varrho_{d}^0(t,x-z)+[\hat q_\alpha(s,z,y)]^2\varrho_{d}^0(s,z-y)\Big),
\end{align*}
where in the last inequality we have used the fact
\begin{align}
\hat q_\alpha(t,z,x)\preceq \hat q_\alpha(s,z,y)\,\,\Leftrightarrow\,\,\varrho_{d}^0(t,x-z)\preceq \varrho_{d}^0(s,z-y) \label{semi}
\end{align}
and the symmetry in $x$ and $y$. The proof is finished.
\end{proof}

\subsection{Gradient estimates}
In this subsection, we derive gradient estimates for $r^D(t,x,y)$.
Recall that $r^D(t,x,y)$ is the transition density of $Y^D$. By the construction of $Y^D$,
it holds (see \cite[(2.2)]{Song}) that
\begin{align}
r^D(t,x,y)=\int_0^\infty p^D_{2}(s,x,y)\mu(t,s)\dif s, \label{rtxy}
\end{align}
where $p^{D}_{2}(t,x,y)$ is the Dirichlet heat kernel of $\Delta|_{D}$, and $\mu(t,s)$ is the density of the subordinator $T_t$. To derive gradient estimates for $r^D(t,x,y)$, we need to recall some estimates for $p^D_2(t,x,y)$.

For any $\gamma, \lambda\in \mR$ and $(t, x)\in (0, \infty)\times\mR^d$, we define
$$
\xi^\gamma_{\lambda}(t,x):=t^{-(d+\gamma)/2}\e^{-\lambda |x|^2/t}.
$$
It is known (see \cite[Theorems 3.1 and 3.2]{Song} for instance) there exist constants $\lambda_1,\lambda_2>0$, $C_1>1$ and $C_2<1$ such that
\begin{align}
p^{D}_{2}(t,x,y)&\leq C_1
\left(1\wedge\frac{\rho(x)\rho(y)}{t}\right)\xi^0_{\lambda_1}(t,x-y), \quad (t, x, y)\in (0, \infty)\times D\times D,\label{pd1ub}\\
p^{D}_{2}(t,x,y)&\geq C_2
\left(1\wedge\frac{\rho(x)\rho(y)}{t}\right)\xi^0_{\lambda_2}(t,x-y), \quad  (t, x, y)\in (0, 1]\times D\times D.\label{pd1lb}
\end{align}
Moreover, it follows from \cite[Theorem 2.1]{ZQS} that, for any $T>0$,
there exists a constant $C_T>0$ such that for all $t\in(0,T]$ and $x,y\in D$,
\begin{equation}
|\nabla_x p^D_2(t,x,y)|\leq \left\{\begin{aligned}
&\frac{C_T}{\rho(x)}p^D_2(t,x,y),\qquad\qquad\qquad\quad\quad \textrm{if}\,\,\,\rho(x)\leq \sqrt{t};\\
&\frac{C_T}{\sqrt{t}}\left(1+\frac{|x-y|}{\sqrt{t}}\right)p^D_2(t,x,y),\qquad \textrm{if}\,\,\,\rho(x)>\sqrt{t}.\label{pd3}
\end{aligned}
\right.
\end{equation}
It turns out that  (\ref{pd1ub}), \eqref{pd1lb} and (\ref{pd3}) are not very convenient to use.
To get easy-to-use forms of the estimates above, we first do some manipulations on $p^D_2(t,x,y)$.
We want to separate the terms $\rho(x)$ and $\rho(y)$.
The following elementary observation will be important.

\bl
For any $\lambda_2>\lambda_1>0$ and $\gamma\in\mR$, it holds for all $t>0$ and $x,y\in D$ that
\begin{align}
\left(1\wedge\frac{\rho(x)\rho(y)}{t}\right)\xi^\gamma_{\lambda_2}(t,x-y)\preceq \left(1\wedge\frac{\rho(x)}{\sqrt{t}}\right)
\left(1\wedge\frac{\rho(y)}{\sqrt{t}}\right)\xi^\gamma_{\lambda_1}(t,x-y).  \label{pdd}
\end{align}
\el
\begin{proof}
In light of \eqref{elmentary}, it suffices to show that  for any $\lambda_0>0$
\begin{align*}
\Big(\rho(x)+\sqrt{t}\Big)\left(\rho(y)+\sqrt{t}\right)\preceq\big(\rho(x)\rho(y)+t\big)\e^{\lambda_0\frac{|x-y|^{2}}{t}}.
\end{align*}
In fact, using symmetry and the elementary inequality
$$
\rho(x)\leq \rho(y)+|x-y|,
$$
we have
$$
\rho(x)^{2}+\rho(y)^{2}\preceq \rho(x)\rho(y)+|x-y|^{2}.
$$
Thus, we can deduce that
\begin{align*}
\left(\rho(x)+\sqrt{t}\right)\left(\rho(y)+\sqrt{t}\right)&\preceq \rho(x)\rho(y)+t+\rho(x)^{2}+\rho(y)^{2}\\
&\preceq \rho(x)\rho(y)+t+|x-y|^{2}.
\end{align*}
Note that for any  $\lambda_0>0$, we have
$$
|x-y|^{2}\preceq t\cdot\e^{\lambda_0\frac{|x-y|^{2}}{t}}.
$$
The desired result follows immediately.
\end{proof}

Recall the definition of $q_\alpha(t,x,y)$ in (\ref{q}). We give a better form of (\ref{pd1ub})
and \eqref{pd1lb} as follows.
\bl
There exist constants $\lambda_1,\lambda_2>0$, $C_1>1$ and $C_2<1$ such that
\begin{align}
p^{D}_{2}(t,x,y)&\leq C_1q_2(t,x,y)\xi^0_{\lambda_1}(t,x-y), \quad  (t, x, y)\in (0, \infty)\times D\times D,\label{pd2ub}\\
p^{D}_{2}(t,x,y)&\geq C_2q_2(t,x,y)\xi^0_{\lambda_2}(t,x-y), \quad (t,x,y)\in (0, 1]\times D\times D.
\label{pd2lb}
\end{align}
\el
\begin{proof}
The lower bound \eqref{pd2lb} is obvious, we only need to prove the upper bound \eqref{pd2ub}. Combining (\ref{pd1ub})
and \eqref{pd1lb} with (\ref{pdd}), we have that for any $\lambda_0>0$,
\begin{align*}
p^{D}_{2}(t,x,y)&\preceq
\left(1\wedge\frac{\rho(x)}{\sqrt{t}}\right)\left(1\wedge\frac{\rho(y)}{\sqrt{t}}\right)\xi^0_{\lambda_0}(t,x-y).
\end{align*}
Thus, (\ref{pd2ub}) is true when $|x-y|\leq \sqrt{t}$. On the other hand, notice that for $0<\tilde{\lambda}_0<\lambda_0$ we have
\begin{align}\label{new}
\frac{\rho(x)}{\sqrt{t}}\e^{-\lambda_0\frac{|x-y|^{2}}{t}}=\frac{\rho(x)}{|x-y|}\cdot\frac{|x-y|}{\sqrt{t}}\e^{-\lambda_0\frac{|x-y|^{2}}{t}}
\preceq \frac{\rho(x)}{|x-y|}\e^{-\tilde{\lambda}_0\frac{|x-y|^{2}}{t}}.
\end{align}
Combining \eqref{new} with \eqref{pd1ub} gives the desired result for $|x-y|>\sqrt{t}$.
\end{proof}

Now we prove the first and second order gradient estimates for
$p^{D}_{2}(t,x,y)$.

\bl
Let $T>0$. There exist constants $C_T, \lambda_3>0$ such that for $j=1,2$, \\
i) for all $t\in(0,T]$ and $x,y\in D$,
\begin{align}\label{p1}
\left|\nabla^j_x p^{D}_{2}(t,x,y)\right|\leq
C_T\hat q_2(t,y,x)\xi^j_{\lambda_3}(t,x-y);
\end{align}
ii) for all $t\in(T,\infty)$ and $x,y\in D$,
\begin{align}\label{p2}
\left|\nabla^j_x p^{D}_{2}(t,x,y)\right|\leq\frac{C_T}{T^{j/2}}\hat q_2(t,y,x)\xi^0_{\lambda_3}(t,x-y),
\end{align}
where $\nabla^j_x$ denotes the $j$-order derivative with respect to the $x$ variable.
\el
\begin{proof}
For (\ref{p1}), we only need to show that there exist $\lambda_3>0$ and $C_T>0$ such that for every $t\in(0,T]$ and $x, y\in D$,
$$
|\nabla^j_xp^{D}_{2}(t,x,y)|\leq C_T\left(1\wedge\frac{\rho(y)}{\sqrt{t}}\right)\xi^j_{\lambda_3}(t,x-y).
$$
Then applying \eqref{new}, we can get (\ref{p1}).
By \cite[VI.2, Theorem 2.1]{G-M}, we have that for every $t\in(0,T]$ and $x, y\in D$,
$$
|\nabla^j_xp^{D}_{2}(t,x,y)|\preceq \xi^j_{\lambda_3}(t,x-y).
$$
Using the Chapman-Kolmogorov equation, we have
\begin{align*}
|\nabla^j_xp^{D}_{2}(t,x,y)|&\leq \int_D\big|\nabla^j_xp^{D}_{2}(t/2,x,z)\big|\cdot p^{D}_{2}(t/2,z,y)\dif z\\
&\preceq \left(1\wedge\frac{\rho(y)}{\sqrt{t}}\right)\int_D\xi^j_{\lambda_3}(t/2,x-z)\xi^0_{\lambda_1}(t/2,z-y)\dif z\\
&\preceq \left(1\wedge\frac{\rho(y)}{\sqrt{t}}\right)\xi^j_{\lambda_3}(t,x-y).
\end{align*}
Thus (\ref{p1}) is valid. We now prove (\ref{p2}). Similarly, it suffices to show that for every $t>T$ and $x, y\in D$,
$$
\left|\nabla^j_x p^{D}_{2}(t,x,y)\right|\leq
\frac{C_T}{T^{j/2}}\left(1\wedge\frac{\rho(y)}{\sqrt{t}}\right)\xi^0_{\lambda_2}(t,x-y).
$$
By (\ref{pd2ub}), (\ref{p1}) and the Chapman-Kolmogorov equation, we have for $t>T$,
\begin{align}
\left|\nabla^j_x p^{D}_{2}(t,x,y)\right|&\leq\int_{D}\left|\nabla^j_x p^{D}_{2}(T,x,z)\right|\cdot p^{D}_{2}(t-T,z,y)\dif z \no\\
&\leq C_T\int_{\mR^{d}}\xi^j_{\lambda_2}(T,x-z)\xi^0_{\lambda_1}(t-T,z-y)\dif z\no\\
&\leq \frac{C_T}{T^{j/2}}\xi^0_{\lambda_2}(t,x-y).  \label{gr0}
\end{align}
Furthermore, for $t\in(T,2T]$ the same argument yields that
\begin{align*}
\left|\nabla^j_x p^{D}_{2}(t,x,y)\right|&\leq\int_{D}\left|\nabla^j_x p^{D}_{2}(t/2,x,z)\right|\cdot p^{D}_{2}(t/2,z,y)\dif z\\
&\leq C_T\int_{\mR^d}\xi^j_{\lambda_2}(t/2,x-z)\frac{\rho(y)}{\sqrt{t/2}}\xi^0_{\lambda_1}(t/2,z-y)\dif z\\
&\leq \frac{C_T}{T^{j/2}}\rho(y)\xi^1_{\lambda_2}(t,x-y)=\frac{C_T}{T^{j/2}}\frac{\rho(y)}{\sqrt{t}}\xi^0_{\lambda_2}(t, x-y).
\end{align*}
Using (\ref{gr0}) we get that for any $t\in(2T,\infty)$,
\begin{align*}
\left|\nabla^j_x p^{D}_{2}(t,x,y)\right|&\leq\int_{D}\left|\nabla^j_x p^{D}_{2}(t/2,x,z)\right|\cdot p^{D}_{2}(t/2,z,y)\dif z\\
&\leq \frac{C_T}{T^{j/2}}\int_{\mR^d}\xi^0_{\lambda_2}(t/2,x-z)\frac{\rho(y)}{\sqrt{t/2}}\xi^0_{\lambda_1}(t/2,z-y)\dif z\\
&\leq \frac{C_T}{T^{j/2}}\rho(y)\xi^1_{\lambda_2}(t,x-y)=\frac{C_T}{T^{j/2}}\frac{\rho(y)}{\sqrt{t}}\xi^0_{\lambda_2}(t, x-y).
\end{align*}
Combining the above computations, we get the
desired result.
\end{proof}

\br
In fact, in the form of (\ref{pd3}), our result means that for every $t\in(0,T]$,
$$
|\nabla^j_x p^D_{2}(t,x,y)|\leq
C_T\frac{(|x-y|+\sqrt{t})^{1-j}}{\rho(x)\wedge(|x-y|+\sqrt{t})}q_2(t,x,y)\xi^0_{\lambda_3}(t,x-y).
$$
Compared with (\ref{pd1ub}), (\ref{pd1lb}) and (\ref{pd3}),
the additional term $|x-y|$ in (\ref{pd2ub})--(\ref{pd2lb}) and (\ref{p1})--(\ref{p2}) is of critical importance
in our derivation of the gradient estimates of $r^D(t,x,y)$ below.
\er

Recall the definition of $q^D(t,x,y)$ in (\ref{qq}). Now, we are ready to derive the following gradient estimates for the Dirichlet heat kernel
$r^D(t,x,y)$.

\bl
Let $T>0$. There exists a constant $C_T>0$ such that for $j=1,2$,
all $t\in(0,T]$ and $x,y\in D$,
\begin{align}\label{grad}
|\nabla^j_x
r^D(t,x,y)|\leq C_T\frac{(|x-y|+t^{1/\alpha})^{1-j}}{\rho(x)\wedge(|x-y|+t^{1/\alpha})}q^D(t,x,y).
\end{align}
Moreover, for any $\vartheta\in(0,1)$ and $t\in(0,T]$, $x,x',y\in D$, we have
\begin{align}
|\nabla_x r^D(t,x,y)-\nabla_x r^D(t,x',y)| \leq C_T|x-x'|^\vartheta \hat q_{\alpha}(t,y,\widetilde x)\varrho_{d+1+\vartheta}^1(t,\widetilde x-y), \label{hold}
\end{align}
where $\widetilde x$ is the point among $x$ and $x'$ which is closer to $y$.
\el
\begin{proof}
We claim that for $j=1,2$,
\begin{align}
|\nabla^j_x r^D(t,x,y)|\preceq
\hat q_\alpha(t,y,x)\varrho^1_{d+j}(t,x-y).\label{gr22}
\end{align}
As a consequence of this claim, we get
\begin{align*}
|\nabla^j_x r^D(t,x,y)|&\preceq\frac{1}{(|x-y|+t^{1/\alpha})^j\hat q_\alpha(t,x,y)} \hat q_\alpha(t,x,y)\hat q_\alpha(t,y,x)\varrho_{d}^1(t,x-y)\\
&\asymp\frac{(|x-y|+t^{1/\alpha})^{1-j}}{\rho(x)\wedge(|x-y|+t^{1/\alpha})}q^D(t,x,y).
\end{align*}
Now we prove the claim (\ref{gr22}).
From \cite[(4.1)]{Song}, we know that for all
$\xi\in\mathbb{R}^{d},$
\begin{align*}
\int_{0}^{\infty}\!s^{-d/2}\e^{-\frac{|\xi|^{2}}{s}}\mu(t,s)\dif
s\asymp\varrho^1_{d}(t,\xi).
\end{align*}
Combining this with (\ref{rtxy}), (\ref{p1}) and (\ref{p2}), we can
get
\begin{align}
|\nabla^j_x r^D(t,x,y)|&\leq \int_{0}^{1}\left|\nabla^j_x p^{D}_{2}(s,x,y)\right|\mu(t,s)\dif s+\int_{1}^{\infty}\left|\nabla^j_x p^{D}_{2}(s,x,y)\right|\mu(t,s)\dif s\no\\
&\preceq\left(1\wedge\frac{\rho(y)}{|x-y|}\right)\left[\int_{0}^{\infty}\!\xi^j_{\lambda_3}(s,x-y)\mu(t,s)\dif s
+\int_{0}^{\infty}\!\xi^0_{\lambda_3}(s,x-y)\mu(t,s)\dif s\right]\no\\
&\asymp\left(1\wedge\frac{\rho(y)}{|x-y|}\right)\Big[\varrho^1_{d+j}(t,x-y)+\varrho_{d}^1(t,x-y)\Big]\no\\
&\preceq\left(1\wedge\frac{\rho(y)}{|x-y|}\right)\varrho^1_{d+j}(t,x-y),
\label{e1}
\end{align}
where in the last inequality we have used the fact that $D$ is bounded and $t\in(0,T]$. Thus, (\ref{gr22}) is true when $|x-y|\geq t^{1/\alpha}$. For the case that $|x-y|<t^{1/\alpha}$,
we may argue similarly to get that
\begin{align*}
|\nabla^j_x r^D(t,x,y)|
&\preceq \rho(y)\left[\int_{0}^{\infty}\!\xi^{j+1}_{\lambda_3}(s,x-y)\mu(t,s)\dif s
+\int_{0}^{\infty}\!\xi^1_{\lambda_3}(s,x-y)\mu(t,s)\dif s\right]\\
&\asymp\rho(y)\Big[\varrho^1_{d+j+1}(t,x-y)+\varrho_{d+1}^1(t,x-y)\Big]\\
&\preceq\frac{\rho(y)}{t^{1/\alpha}}\varrho^1_{d+j}(t,x-y).
\end{align*}
This together with estimate (\ref{e1}) implies (\ref{gr22}).

For (\ref{hold}), without loss of generality, we may assume that $|x-y|\leq |x'-y|$. Using (\ref{gr22}) with $j=1$, we can get that when $|x-x'|\geq (|x-y|+t^{1/\alpha})/2$,
\begin{align*}
\sQ&:=|\nabla_x r^D(t,x,y)-\nabla_x r^D(t,x',y)|\\
&\leq C_T|x-x'|^\vartheta(|x-y|+t^{1/\alpha})^{-\vartheta}\Big(\hat q_\alpha(t,y,x)\varrho^{1}_{d+1}(t,x-y)\\
&\qquad\qquad\qquad\qquad\qquad\quad\qquad\quad+\hat q_\alpha(t,y,x')\varrho^{1}_{d+1}(t,x'-y)\Big)\\
&\leq C_T|x-x'|^\vartheta\hat q_\alpha(t,y,x)\varrho_{d+1+\vartheta}^{1}(t,x-y).
\end{align*}
When $|x-x'|<(|x-y|+t^{1/\alpha})/2$, we have by the mean value theorem and (\ref{gr22}) with $j=2$ that for some $\varepsilon\in[0,1]$,
\begin{align*}
\sQ&\leq C_T|x-x'|\hat q_\alpha\big(t,y,x+\varepsilon(x'-x)\big)\varrho^{1}_{d+2}\big(t,x+\varepsilon(x'-x)-y\big)\\
&\leq C_T|x-x'|\hat q_\alpha(t,y,x)\varrho_{d+2}^{1}(t, x-y)\\
&\leq C_T|x-x'|^\vartheta\hat q_\alpha(t,y,x)\varrho_{d+1+\vartheta}^{1}(t,x-y).
\end{align*}
The proof is finished.
\end{proof}

\section{Proof of Theorem \ref{main}}

 Let
\begin{align*}
\widehat{\bf K}_{D}^{\alpha-1}:=\bigg\{f\in L^1_{loc}(D): \lim_{t\downarrow0}\sup_{x\in D}\int_{D}&\left(1\wedge \frac{\rho(y)}{|x-y|}\right)\\
&\times\bigg(\frac{1}{|x-y|^{d+1-\alpha}}\wedge\frac{t^2}{|x-y|^{d+\alpha+1}}\bigg)|f(y)|\dif y=0\bigg\}.
\end{align*}
We first give the following result about our Kato class.

\bl\label{include}
We have ${\bf K}_D^{\alpha-1}\subset \widehat{\bf K}_D^{\alpha-1}\subset\mK_D^{0}$.
Moreover, for any $\gamma\geq 0$, if $1<p,q\leq \infty$ satisfy
\begin{align}
\frac{d}{\alpha p}+\frac{1}{q}<1-\frac{1+\gamma}{\alpha},    \label{pq}
\end{align}
then $L^q(\mR;L^p(D))\subseteq\mK_D^\gamma$.
\el
\begin{proof}
It follows from \cite[Lemma 2.1]{C-S-X-X}, which follows from \cite[Corollary 12]{Bo-Ja0},
that a real-valued function $f$ belongs to ${\bf K}^{\alpha-1}_D$ if and only if
\begin{equation*}
\lim_{t\rightarrow0}\sup_{x\in D}\int_{D}\bigg(\frac{1}{|x-y|^{d+1-\alpha}}\wedge\frac{t^2}{|x-y|^{d+\alpha+1}}\bigg)|f(y)|\dif y=0.
\end{equation*}
Thus the first inclusion is obvious.
To show that a real-valued time-independent function $f$ on $D$
belongs to $\mK_D^0$,
it suffices to show that
$$
\int_0^t\varrho_{d+1}^1(s,x-y)\dif s\preceq \frac{1}{|x-y|^{d+1-\alpha}}\wedge\frac{t^2}{|x-y|^{d+\alpha+1}}.
$$
This follows directly from \cite[Lemma 2.3]{C-S-X-X} with $\gamma=1$.
Thus the second inclusion is valid. Now we prove the third inclusion. By H\"older's inequality, we get
$$
K_f^\gamma(\delta)\leq\Bigg(\int_{\mR}\left(\int_{D}|f(s,y)|^p
\dif y\right)^{\frac{q}{p}}\dif s\Bigg)^{\frac{1}{q}}I_{\alpha,\gamma}(\delta),
$$
where
$$
I_{\alpha,\gamma}(\delta):=\delta^{\frac{\gamma}{\alpha}}\Bigg(\int^\delta_0\big[s^{-\gamma/\alpha}+(\delta-s)^{-\gamma/\alpha}\big]^{q^*}\Bigg(\!\int_{\mR^d}\frac{s^{p^*}}
{\big(|y|+s^{1/\alpha}\big)^{(d+\alpha+1)p^{*}}}\dif y\Bigg)^{\frac{q^*}{p^*}}\dif s\Bigg)^{\frac{1}{q^{*}}},
$$
with $q^*:=\frac{q}{q-1}$ and $p^*:=\frac{p}{p-1}$. Noticing that
\begin{align*}
\int_{\mR^d}\frac{s^{p^*}}{\big(|y|+s^{1/\alpha}\big)^{(d+\alpha+1)p^{*}}}\dif y
&\leq s^{p^*}\left(\int_{|y|\leq s^{1/\alpha}}s^{-\frac{(d+\alpha+1)p^{*}}{\alpha}}\dif y
+ \int_{|y|> s^{1/\alpha}}\frac{\dif y}{|y|^{(d+\alpha+1)p^{*}}}\right)\\
&\preceq s^{\frac{d-(d+1)p^{*}}{\alpha}},
\end{align*}
we have
$$
I_{\alpha,\gamma}(\delta)
\preceq\delta^{\frac{\gamma}{\alpha}}\Bigg(\int^\delta_0\big[s^{-\gamma/\alpha}+(\delta-s)^{-\gamma/\alpha}\big]^{q^*}s^{\frac{dq^*}{\alpha p^*}-\frac{(d+1)q^*}{\alpha}}\dif s\Bigg)^{\frac{1}{q^{*}}}.
$$
Thus $I_{\alpha,\gamma}(\delta)$ converges to zero as $\delta\to 0$ provided that
$$
-\frac{\gamma q^*}{\alpha}+\frac{dq^*}{\alpha p^*}-\frac{d+1}{\alpha}q^*+1>0\Leftrightarrow (\ref{pq}).
$$
The desired result follows.
\end{proof}

The following lemma is related to the smallness of $b\cdot\nabla$ as a perturbation of $-(-\Delta|_{D})^{\alpha/2}$, which plays an important role in proving our main result.

\bl\label{integral1}
Let $\delta>0$ and $b\in \mK_{D}^{0}$. Then for all $0\leq s<t\leq s+\delta$ and $x, y\in
D$, we have
\begin{align*}
\int_s^t\!\!\!\int_{D}r^D(r-s,x,z)|b(r,z)|\cdot|\nabla_z
r^D(t-r,z,y)|\dif z\dif r \leq C(\delta)r^D(t-s,x,y),
\end{align*}
where $C(\delta)$ is a positive constant with $C(\delta)\rightarrow0$ as $\delta\downarrow0$.
\el
\begin{proof}
In this proof we always assume that $0\leq s<t\leq s+\delta$ and $x, y\in D$.
For brevity, we write
$$
\cW:=\frac{r^D(r-s,x,z)|\nabla_z r^D(t-r,z,y)|}{r^D(t-s,x,y)}.
$$
It follows from (\ref{grad}) that
\begin{align*}
\cW&\preceq \frac{r^D(r-s,x,z)r^D(t-r,z,y)}{r^D(t-s,x,y)}\cdot\frac{1}{\rho(z)\wedge\left(|z-y|+(t-r)^{1/\alpha}\right)}\\
&\preceq \frac{r^D(r-s,x,z)r^D(t-r,z,y)}{r^D(t-s,x,y)}\left(\frac{1}{\rho(z)}+\frac{1}{|z-y|+(t-r)^{1/\alpha}}\right)=:\cW_1+\cW_2.
\end{align*}
By (\ref{3r}), we have that
\begin{align}\label{i}
\cW_1&\preceq \big((r-s)\wedge (t-r)\big)\Big([\hat q_\alpha(r-s,z,x)]^2\varrho_{d}^0(r-s,x-z)\no\\
&\qquad\qquad\qquad\qquad\qquad+[\hat q_\alpha(t-r,z,y)]^2\varrho_{d}^0(t-r,z-y)\Big)\frac{1}{\rho(z)}\no\\
&\preceq \hat q_\alpha(r-s,z,x)\varrho^1_{d+1}(r-s,x-z)
+\hat q_\alpha(t-r,z,y)\varrho^1_{d+1}(t-r,z-y).
\end{align}
Again by (\ref{3r}), we have
\begin{align*}
\cW_2&\preceq \big((r-s)\wedge (t-r)\big)\Big([\hat q_\alpha(r-s,z,x)]^2\varrho_{d}^0(r-s,x-z)\\
&\qquad\qquad\qquad\qquad+[\hat q_\alpha(t-r,z,y)]^2\varrho_{d}^0(t-r,z-y)\Big)\frac{1}{|z-y|+(t-r)^{1/\alpha}}.
\end{align*}
By the same argument as in (\ref{semi}), in the case $|x-z|+(r-s)^{1/\alpha}\leq |z-y|+(t-r)^{1/\alpha}$, we have
\begin{align*}
\cW_2&\preceq [\hat q_\alpha(r-s,z,x)]^2\varrho_{d}^1(r-s,x-z)\frac{1}{|x-z|+(r-s)^{1/\alpha}}\\
&\preceq \hat q_\alpha(r-s,z,x)\varrho^1_{d+1}(r-s,x-z).
\end{align*}
In the case $|x-z|+(r-s)^{1/\alpha}>|z-y|+(t-r)^{1/\alpha}$, we have
\begin{align*}
\cW_2&\preceq  [\hat q_\alpha(t-r,z,y)]^2\varrho_{d}^1(t-r,z-y)\frac{1}{|z-y|+(t-r)^{1/\alpha}}\\
&\preceq \hat q_\alpha(t-r,z,y)\varrho^1_{d+1}(t-r,z-y).
\end{align*}
Hence,
$$
\cW_2\preceq
\hat q_\alpha(r-s,z,x)\varrho^1_{d+1}(r-s,x-z)+\hat q_\alpha(t-r,z,y)\varrho^1_{d+1}(t-r,z-y),
$$
which together with (\ref{i}) yields that
$$
\cW\preceq \hat q_\alpha(r-s,z,x)\varrho^1_{d+1}(r-s,x-z)
+\hat q_\alpha(t-r,z,y)\varrho^1_{d+1}(t-r,z-y).
$$
Consequently, by the definition of Kato class
$\mK_{D}^{0},$ it holds that
\begin{align*}
&\quad\int_s^t\!\!\!\int_{D}r^D(r-s,x,z)|b(r,z)|\cdot|\nabla_{z}
r^D(t-r,z,y)|\dif z\dif r\\
&\preceq\int_s^t\!\!\!\int_{D}\hat q_\alpha(r-s,z,x)\varrho^1_{d+1}(r-s,x-z)|b(r,z)|\dif
z\dif r\cdot r^D(t-s,x,y)\\
&\quad+\int_s^t\!\!\!\int_{D}\hat q_\alpha(t-r,z,y)\varrho^1_{d+1}(t-r,z-y)|b(r,z)|\dif
z\dif r\cdot r^D(t-s,x,y)\\
&\leq 2K_{b}^{0}(\delta)r^D(t-s,x,y),
\end{align*}
where $K_{b}^{0}(\delta)$
is defined in Definition \ref{Definition}. The proof is thus finished.
\end{proof}

To derive the gradient estimate of the Dirichlet heat kernel, we shall also need the following
result.

\bl\label{integral2}
Let $\delta>0$ and $b\in \mK_{D}^{0}$. Then for all $0\leq s<t\leq s+\delta$ and $x, y\in D$, we have
\begin{align*}
&\int_s^t\!\!\!\int_{D}|\nabla_x r^D(r-s,x,z)||b(r,z)|\cdot|\nabla_zr^D(t-r,z,y)\dif z\dif r\\
&\leq
\hat C(\delta)\frac{1}{\rho(x)\wedge(|x-y|+(t-s)^{1/\alpha})}r^D(t-s,x,y),
\end{align*}
where $\hat C(\delta)$ is a positive constant with $\hat C(\delta)\rightarrow0$ as $\delta\downarrow 0$.
\el
\begin{proof}
In this proof we always assume that $0\leq s<t\leq s+\delta$ and $x, y\in D$.
Define
\begin{align*}
\cV&:=\frac{|\nabla_{x}r^D(r-s,x,z)|\cdot|\nabla_{z}
r^D(t-r,z,y)|}{\hat q_\alpha(t-s,y,x)\varrho_{d+1}^1(t-s,x-y)}.
\end{align*}
It follows from (\ref{gr22}) that
$$
\cV\preceq \cQ\cdot\frac{\varrho_{d+1}^1(r-s,x-z)\varrho_{d+1}^1(t-r,z-y)}{\varrho_{d+1}^1(t-s,x-y)},
$$
where
\begin{align*}
\cQ&:=\frac{\hat q_\alpha(r-s,z,x)\hat q_\alpha(t-r,y,z)}{\hat q_\alpha(t-s,y,x)}.
\end{align*}
Using \eqref{new}, we get
\begin{align*}
\cQ&\asymp\rho(z)\cdot\frac{\rho(y)+|x-y|+(t-s)^{1/\alpha}}{(\rho(z)+|x-z|+(r-s)^{1/\alpha})(\rho(y)+|z-y|+(t-r)^{1/\alpha})}\\
&\preceq\rho(z)\cdot\frac{\rho(z)+|x-z|+(r-s)^{1/\alpha}+\rho(z)+|z-y|+(t-r)^{1/\alpha}}{(\rho(z)+|x-z|+(r-s)^{1/\alpha})(\rho(z)+|z-y|+(t-r)^{1/\alpha})}\\
&=\frac{\rho(z)}{\rho(z)+|x-z|+(r-s)^{1/\alpha}}+\frac{\rho(z)}{\rho(z)+|z-y|+(t-r)^{1/\alpha}}\\
&\asymp
\Big(\hat q_{\alpha}(r-s,z,x)+\hat q_{\alpha}(t-r,z,y)\Big).
\end{align*}
Combining this with (\ref{333}), and by the same argument as in (\ref{semi}), we further have that
\begin{align}
\cV&\preceq[(r-s)\wedge (t-r)]\Big(\hat q_{\alpha}(r-s,z,x)+\hat q_{\alpha}(t-r,z,y)\Big)\no\\
&\qquad\qquad\qquad\qquad\qquad\qquad\times \big(\varrho^0_{d+1}(r-s,x-z)+\varrho^0_{d+1}(t-r,z-y)\big)\no\\
&\preceq \hat q_\alpha(r-s,z,x)\varrho^1_{d+1}(r-s,x-z)
+\hat q_\alpha(t-r,z,y)\varrho^1_{d+1}(t-r,z-y).   \label{3g}
\end{align}
Hence,
\begin{align*}
&\int_s^t\!\!\!\int_{D}|\nabla_x r^D(r-s,x,z)||b(r,z)|\cdot|\nabla_zr^D(t-r,z,y)\dif z\dif r\\
&\preceq
K_b^0(\delta)\hat q_\alpha(t-s,y,x)\varrho_{d+1}^1(t-s,x-y)\\
&\leq K_b^0(\delta)\frac{1}{\rho(x)\wedge(|x-y|+(t-s)^{1/\alpha})}r^D(t-s,x,y),
\end{align*}
which yields the desired result. The proof is finished.
\end{proof}

We now proceed to solve the integral equation (\ref{Duhamel}). For
all $0\leq s<t$ and $x,y\in D$, set $r_0(s,x;t,y):=r^D(t-s,x,y)$,
and define inductively that for $k\geq1$,
\begin{align}\label{induction}
r_k(s,x;t,y):=\int_s^t\!\!\!\int_Dr_{k-1}(s,x;r,z)b(r,z)\cdot\nabla_z
r_{0}(r,z;t,y)\dif z\dif r.
\end{align}
The following result is an easy consequence of Lemma \ref{integral1}
and Lemma \ref{integral2}.

\bl
Let $\delta>0$ and $b\in\mK^{0}_D$.
Then there exits a constant $c_1>1$ such that for all $k\geq 1$,
$0\leq s<t\leq s+\delta$ and $x,y\in D$, we have
\begin{align}\label{3.5}
|r_k(s,x;t,y)|\leq [c_1C(\delta)]^{k}q^D(t-s,x,y)
\end{align}
and
\begin{align}\label{GK}
|\nabla_x r_k(s,x;t,y)|\leq [c_1\hat C(\delta)]^{k}\frac{1}{\rho(x)\wedge(|x-y|+(t-s)^{1/\alpha})}q^D(t-s,x,y),
\end{align}
where $C(\delta)$ and $\hat C(\delta)$ are the constants in Lemma
\ref{integral1} and Lemma \ref{integral2}, respectively. Moreover,
it holds that
\begin{align}
r_k(s,x;t,y)=\int_s^t\!\!\!\int_Dr_{0}(s,x;r,z)b(r,z)\cdot\nabla_z
r_{k-1}(r,z;t,y)\dif z\dif r.    \label{pkpk}
\end{align}
\el
\begin{proof}
We first prove (\ref{3.5}) by induction. By Lemma
\ref{integral1} and the definition of $q^D(t,x,y)$, we know that (\ref{3.5}) holds for $k=1$. Now
suppose that it holds for  $k>1$. Then by definition and
using Lemmas \ref{integral1} and \ref{sharp}, we have
\begin{align*}
|r_{k+1}(s,x;t,y)|
&\leq\int_s^t\!\!\!\int_D|r_{k}(s,x;r,z)|\cdot|b(r,z)|\cdot|\nabla_z
r_{0}(r,z;t,y)|\dif z\dif r\\
&\leq
[c_1C(\delta)]^{k}\int_s^t\!\!\!\int_Dr^{D}(s,x;r,z)|b(r,z)|\cdot|\nabla_z
r_{0}(r,z;s,y)|\dif z\dif r\\
&\leq [c_1C(\delta)]^{k+1}q^D(t-s,x,y).
\end{align*}
Following the same argument with Lemma \ref{integral1} replaced by
Lemma \ref{integral2}, we can show that (\ref{GK}) is true. We
proceed to prove (\ref{pkpk}). It is obvious that (\ref{pkpk}) holds
for $k=1$. Suppose that it is true for $k>1$. Then, we have
by (\ref{induction}) and Fubini's theorem that
\begin{align*}
r_{k+1}(s,x;t,y)&=\int_s^t\!\!\!\int_Dr_{k}(s,x;r,z)b(r,z)\cdot\nabla_z r_{0}(r,z;s,y)\dif z\dif r\\
&=\int_s^t\!\!\!\int_D\int_s^r\!\!\!\int_Dr_{0}(s,x;r',z')b(r',z')\cdot\nabla_{z'} r_{k-1}(r',z';r,z)\dif z'\dif r'\\
&\qquad\qquad\qquad\qquad\qquad\qquad\qquad\times b(r,z)\cdot\nabla_z r_{0}(r,z;t,y)\dif z\dif r\\
&=\int_s^t\!\!\!\int_Dr_{0}(s,x;r',z')b(r',z')\cdot\!\int_{r'}^{t}\!\!\int_D\nabla_{z'} r_{k-1}(r',z';r,z)\\
&\qquad\qquad\qquad\qquad\qquad\qquad\qquad\times b(r,z)\cdot\nabla_z r_{0}(r,z;t,y)\dif z\dif r\dif z'\dif r'\\
&=\int_s^t\!\!\!\int_Dr_{0}(s,x;r',z')b(r',z')\cdot\nabla_{z'}r_{k}(r',z';t,y)\dif z'\dif r'.
\end{align*}
The proof is complete.
\end{proof}

Now, we are ready to give:

\begin{proof}[Proof of Theorem \ref{main}]
Let $r_k$ be defined by \eqref{induction}. For $\delta>0$, define
$\mD_\delta:=\{(s,x;t,y): x,y\in D\,\, \text{and}\,\,0\leq s<t\leq s+\delta\}$.
It follows from Lemma \ref{integral1} that there exists a
$\delta_0\in (0, 1]$ such that for all $0\leq s<t\leq s+\delta_0$, we have
$c_1C(\delta_0)<1/4$,
where $c_1$ and $C(\delta_0)$ are the constants from Lemma \ref{induction}. Hence,
\begin{align}
\sum_{k=0}^{\infty}|r_{k}(s,x;t,y)|\leq
\frac43 q^{D}(t-s,x,y)\,\,\,\text{on}\,\,\,\mD_{\delta_0},   \label{es}
\end{align}
which means that the series $\sum_{k=0}^{\infty}r_{k}(s,x;t,y)$
converges on $\mD_{\delta_0}$. Define
$r^{D,b}(s,x;t,y):=\sum_{k=0}^{\infty}r_{k}(s,x;t,y)$ on $\mD_{\delta_0}$.
By \eqref{induction}, we have
\begin{equation*}
\sum_{k=0}^{n+1}r_k(s,x;t,y)=r_0(s,x;t,y)+\int_s^t\!\!\!\int_D\sum_{k=0}^nr_k(s,x;r,z)b(r,z)\cdot\nabla_z r_0(r,z;t,y)\dif z\dif r.
\end{equation*}
Letting $n\rightarrow\infty$ on both sides, we get (\ref{Duhamel}).

\vspace{2mm} \noindent(i) The upper bound on $\mD_{\delta_0}$ follows by (\ref{es}). As
for the lower bound on $\mD_{\delta_0}$, we have
\begin{align*}
r^{D,b}(s,x;t,y)\geq r^D(t-s,x,y)-\sum_{k=1}^\infty |r_k(s,x;t,y)|\geq
\frac{2}{3}r^D(t-s,x,y).
\end{align*}
Thus, (\ref{estimate}) is valid on $\mD_{\delta_0}$.

Now let $\tilde{r}^{D,b}(s,x;t,y)$ be another solution to
the integral equation (\ref{Duhamel}) satisfying (\ref{estimate}) on $\mD_{\delta_0}$.
We claim that for every $k\in \mathbb{N}$, there exists a constant $C_0$ such that
on $\mD_{\delta_0}$,
\begin{align}\label{unique}
|r^{D,b}(s,x;t,y)-\tilde{r}^{D,b}(s,x;t,y)|\leq C_0[c_1C(\delta_0)]^{k}q^{D}(t-s,x,y).
\end{align}
Indeed, for $k=1$, using (\ref{Duhamel}), (\ref{estimate}) and Lemma
\ref{integral1} we have
\begin{align*}
&\quad|r^{D,b}(s,x;t,y)-\tilde{r}^{D,b}(s,x;t,y)|\\
&\leq\int_{s}^{t}\!\!\!\int_{D}\big(|r^{D,b}(s,x;r,z)|+|\tilde{r}^{D,b}(s,x;r,z)|\big)\cdot|b(r,z)|\cdot|\nabla_z
r^{D}(t-r,z,y)|\dif z\dif r\\
&\leq C_0\int_{s}^{t}\!\!\!\int_{D}r^{D}(r-s,x,z)\cdot|b(r,z)|\cdot|\nabla_z
r^{D}(t-r,z,y)|\dif z\dif r\leq C_0c_1C(\delta_0)q^{D}(t-s,x,y).
\end{align*}
Suppose that \eqref{unique} holds for some $k\in \mathbb{N}$. By
(\ref{Duhamel}), Lemma \ref{integral1} and the induction hypothesis,
we have
\begin{align*}
&\quad|r^{D,b}(s,x;t,y)-\tilde{r}^{D,b}(s,x;t,y)|\\
&\leq\int_{s}^{t}\!\!\!\int_{D}|r^{D,b}(s,x;r,z)-\tilde{r}^{D,b}(s,x;r,z)|\cdot|b(r,z)|\cdot|\nabla_z
r^D(r,z;t,y)|\dif z\dif r\\
&\leq
C_0[c_1C(\delta_0)]^{k}\int_{s}^{t}\!\!\!\int_{D}r^D(r-s,x,z)\cdot|b(r,z)|\cdot|\nabla_z
r^D(t-r,z,y)|\dif z\dif r\\
&\leq C_0[c_1C(\delta_0)]^{k+1}q^D(t-s,x,y).
\end{align*}
Since $c_1C(\delta_0)<1$, letting $k\rightarrow\infty$ we obtain the
uniqueness.

\vspace{2mm} \noindent(ii)
By choosing $\delta_0$ smaller if necessary, we can assume that $c_1\hat C(\delta_0)<1$
for $0\leq s<t\leq s+\delta_0$,
where $c_1$ and $\widehat{C}(\delta_0)$ are the constants from Lemma \ref{induction}.
It then follows from \eqref{GK} that on $\mD_{\delta_0}$,
\begin{align*}
\left|\sum_{k=0}^{\infty}\nabla_x
r_{k}(s,x;t,y)\right|\preceq\frac{1}{\rho(x)\wedge(|x-y|+(t-s)^{1/\alpha})}q^D(t-s,x,y),
\end{align*}
which means that (\ref{vi}) is true. Moreover, by (\ref{pkpk}) and
Fubini's theorem, we have
\begin{align*}
&\quad r^{D,b}(s,x;t,y)=\sum_{k=0}^{\infty}r_{k}(s,x;t,y)\\
&=r^D(s,x;t,y)+\sum_{k=0}^{\infty}\int_s^t\!\!\!\int_Dr_{0}(s,x;r,z)b(r,z)\cdot\nabla_z
r_{k}(r,z;t,y)\dif z\dif r\\
&=r^D(s,x;t,y)+\int_s^t\!\!\!\int_Dr_{0}(s,x;r,z)b(r,z)\cdot\nabla_z
r^{D,b}(r,z;t,y)\dif z\dif r.
\end{align*}
This yields (\ref{duhamel}).

\vspace{2mm} \noindent(iii) By Fubini's theorem, we have
$$
\int_Dr^{D,b}(s,x;r,z)r^{D,b}(r,z;t,y)\dif
z=\sum_{n=0}^{\infty}\sum_{m=0}^n\int_Dr_m(s,x;r,z)r_{n-m}(r,z;t,y)\dif
z.
$$
Thus, for proving (\ref{eq21}), it suffices to show that for each
$n\in\mN_0$,
\begin{align}
\sum_{m=0}^n\int_Dr_{m}(s,x;r,z)r_{n-m}(r,z;t,y)\dif z=r_n(s,x;t,y).
\label{ck}
\end{align}
It is clear that the above equality holds for $n=0$. Suppose now
that it holds for some $n\in\mN$. Write
$$
\sum_{m=0}^{n+1}\int_Dr_{m}(s,x;r,z)r_{n+1-m}(r,z;t,y)\dif
z=\cJ_1+\cJ_2,
$$
where
$$
\cJ_1:=\int_Dr_{n+1}(s,x;r,z)p_{0}(r,z;t,y)\dif z
$$
and
$$
\cJ_2:=\sum_{m=0}^{n}\int_Dr_{m}(s,x;r,z)p_{n+1-m}(r,z;t,y)\dif z.
$$
By (\ref{induction}) and Fubini's theorem, we have
\begin{align*}
\cJ_{1}&=\int_{D}\left(\int_{s}^{r}\!\!\!\int_{D}r_{n}(s,x;r',z')b(r',z')\cdot
\nabla_{z'} r_{0}(r',z';r,z)\dif z'\dif r'\right)r_{0}(r,z;t,y)\dif z\\
&=\int_{s}^{r}\!\!\!\int_{D}r_{n}(s,x;r',z')b(r',z')\cdot\left(\int_{D}\nabla_{z'}
r_{0}(r',z';r,z)r_{0}(r,z;t,y)\dif z\right)\dif z'\dif r'\\
&=\int_{s}^{r}\!\!\!\int_{D}r_{n}(s,x;r',z')b(r',z')\cdot\nabla_{z'}
r_{0}(r',z';t,y)\dif z'\dif r'.
\end{align*}
Similarly, by (\ref{induction}) and the induction hypothesis, we
have
\begin{align*}
\cJ_{2}=\int_{r}^{t}\!\!\!\int_{D}r_{n}(s,x;r',z')b(r',z')\cdot\nabla_{z'}
r_{0}(r',z';t,y)\dif z'\dif r'.
\end{align*}
Hence,
\begin{align*}
\cJ_{1}+\cJ_{2}=\int_{s}^{t}\!\!\!\int_{D}r_{n}(s,x;r',z')b(r',z')\cdot\nabla_{z'}
r_{0}(r',z';t,y)\dif z'\dif r'=r_{n+1}(s,x;t,y),
\end{align*}
which gives (\ref{ck}).

\vspace{2mm}
Now, we can extend  $r^{D,b}(s,x;t,y)$ from $\mD_{\delta_0}$ to
$\{(s,x;t,y): x,y\in D\,\, \text{and}\,\,0\leq s<t<\infty\}$.
Then it is routine to
extend the above assertions (i)-(iii) on $\mD_{\delta_0}$ to
$D_{\delta}$ for any$\delta>0$.

\vspace{2mm} \noindent(iv) Let $R_{s,t}f(x):=\int_Dr^D(t-s,x,y)f(y)\dif y$.
By (\ref{Duhamel}), we have for any $f\in C_c^2(D)$,
\begin{align}\label{RT5}
R_{s,t}^{D,b}f(x)=R_{s,t}f(x)+\int^t_s R_{s,r}^{D,b}(b\cdot\nabla R_{r,t}f)(x)\dif r.
\end{align}
It then follows that
\begin{align}
R_{s,t}^{D,b}f(x)-f(x)&=R_{s,t}f(x)-f(x)+\int^t_s R_{s,r}^{D,b}(b\cdot\nabla R_{r,t}f)(x)\dif r\no\\
&=\int_s^t R_{s,r}(-(-\Delta|_D)^{\alpha/2})f(x)\dif r+\int^t_s R_{s,r}^{D,b}(b\cdot\nabla R_{r,t}f)(x)\dif r,\label{RT4}
\end{align}
and, by \eqref{RT5} and Fubini's theorem,
\begin{align*}
&\quad\int_{s}^{t}R^{D,b}_{s,r}(-(-\Delta|_D)^{\alpha/2})f(x)\dif
s-\int_s^t R_{s,r}(-(-\Delta|_D)^{\alpha/2})f(x)\dif r\\
&=\int_{s}^{t}\!\!\int_{s}^{r}R^{D,b}_{s,u}\big(b\cdot\nabla R_{u,r}(-(-\Delta|_D)^{\alpha/2})f\big)(x)\dif u\dif r\\
&=\int_{s}^{t}R^{D,b}_{s,u}b\cdot\nabla\left(\int_{u}^{t}R_{u,r}(-(-\Delta|_D)^{\alpha/2})f(x)\dif r\right)\dif u\\
&=\int_{s}^{t}R^{D,b}_{s,u}b\cdot\nabla\Big(R_{u,t}f(x)-f(x)\Big)\dif u.
\end{align*}
Combining this with \eqref{RT4}, we obtain
\begin{align*}
R_{s,t}^{D,b}f(x)-f(x)=\int_s^tR_{s,r}^{D,b}\sL^{D,b}f(x)\dif
r,
\end{align*}
which gives (\ref{eq23}).

\vspace{2mm} \noindent(v) Since $r^D(t,x,y)$ is the transition density
of the process $Y^D$, so we have for any uniformly continuous function
$f(x)$ with compact supports,
\begin{align*}
\lim_{t\downarrow s}\|R_{s,t}f-f\|_\infty=0.
\end{align*}
Meanwhile, by (\ref{estimate}) and Lemma \ref{integral1} we have
\begin{align*}
&\quad\left|\!\int_{D}\!\!\!\left(\int_{s}^{t}\!\!\!\int_{D}r^{D,b}_{\alpha}(s,x;r,z)b(r,z)\cdot\nabla_z
r^D(t-r,z,y)\dif z\dif r\!\right)\!f(y)\dif y\right|\\
&\preceq
\|f\|_{\infty}\int_{D}\!\!\!\left(\int_{s}^{t}\!\!\!\int_{D}r^D(r-s,x,z)|b(r,z)|\cdot|\nabla_z
r^D(t-r,z,y)|\dif z\dif r\!\right)\dif y\\
&\leq C(\delta)\|f\|_{\infty}\int_{D}r^D(t-s,x,y)\dif y\leq
C(\delta)\|f\|_{\infty},
\end{align*}
which yields (\ref{con}) by (\ref{Duhamel}).

\vspace{2mm} \noindent(vi) Set
$$
\Phi(s,x;t,y):=\int_s^t\!\!\!\int_D r^{D}(r-s,x,z)b(r,z)\cdot\nabla_z r^{D,b}(r,z;t,y)\dif z\dif r.
$$
If we further assume that for $\gamma\in(0,\alpha-1)$, $b\in\mK_D^{\gamma}$, then using (\ref{hold}) we have for any $x,x',y\in D$,
\begin{align*}
&|\Phi(s,x;t,y)-\Phi(s,x';t,y)|\preceq |x-x'|^{\gamma}\int_s^t\!\!\!\int_D\hat q_{\alpha}(r-s,z,\tilde x)\varrho_{d+1+\gamma}^1(r-s,\tilde x-z)\\
&\qquad\qquad\times|b(r,z)|\hat q_{\alpha}(t-r,y,z)\varrho_{d+1}^1(t-r,z-y)\dif z\dif r\\
&\preceq |x-x'|^{\gamma}\int_s^t\!\!\!\int_D(r-s)^{-\gamma/\alpha}\hat q_{\alpha}(r-s,z,\tilde x)\varrho_{d+1}^1(r-s,\tilde x-z)\\
&\qquad\qquad\times|b(r,z)|\hat q_{\alpha}(t-r,y,z)\varrho_{d+1}^1(t-r,z-y)\dif z\dif r\\
&\preceq |x-x'|^{\gamma}\hat q_{\alpha}(t-s,y,\tilde x)\varrho_{d+1}^1(t-s,\tilde x-y)\int_s^t\!\!\!\int_D(r-s)^{-\gamma/\alpha}|b(r,z)|\\
&\qquad\qquad\times\Big(\hat q_{\alpha}(r-s,z,\tilde x)\varrho_{d+1}^1(r-s,\tilde x-z)+\hat q_{\alpha}(t-r,z,y)\varrho_{d+1}^1(t-r,z-y)\Big)\dif z\dif r\\
&\preceq |x-x'|^{\gamma}(t-s)^{-\gamma/\alpha}\hat q_{\alpha}(t-s,y,\tilde x)\varrho_{d+1}^1(t-s,\tilde x-y),
\end{align*}
where the third inequality is due to (\ref{3g}), and the last inequality follows from the definition of $\mK_D^\gamma$. Combining this with (\ref{duhamel}) and (\ref{hold}), we get the desired result. The proof is complete.
\end{proof}

\end{document}